\pdfoutput=1 
\documentclass[letterpaper, 10 pt, conference]{ieeeconf}  

\IEEEoverridecommandlockouts                              

\usepackage{graphicx} 
\usepackage{amsmath} 
\usepackage{amssymb}  

\newtheorem{assumption}{Assumption}

\title{\LARGE \bf
Motion and Communication Co-optimization with Path Planning \\
and  Online Channel Estimation$^{\dag}$\thanks{$^{\dag}$Research supported in part by  the NSF  under Grant Numbers CNS-1239225 and NeTS-1321171.}
\vspace{0.1in}}

\author{Usman Ali$^*$, Hong Cai$^{**}$, Yasamin Mostofi$^{**}$ and Yorai Wardi$^*$
\thanks{$^*$School of Electrical and Computer Engineering, Georgia Institute of Technology, Atlanta, GA 30332. Email: usmanali@gatech.edu, ywardi@ece.gatech.edu.}
\thanks{$^{**}$Department of Electrical and Computer Engineering, University of California, Santa Barbara, CA 93106. Email: hcai@ece.ucsb.edu, ymostofi@ece.ucsb.edu.}
}

\begin{document}

\maketitle
\thispagestyle{empty}
\pagestyle{empty}

\begin{abstract}
 This paper considers the problem of optimally balancing motion energy and communication transmission energy of a mobile robot tasked with
transmitting a given number of data bits to a remote station, while navigating to
a pre-specified destination in  a given amount of time.  The problem is cast in the setting of optimal control, where the
robot has to choose its path, acceleration, and transmission rate along the path so as to minimize its energy required for
transmission and motion, while satisfying various power and communication constraints.
We use realistic models for the robot's channel estimation, motion dynamics, and power  and energy costs. The main contribution of the paper is to show how to co-optimize robot's path along with other communication and motion variables.  Two versions of the
problem are solved: the first is defined offline by assuming that all the channel
measurements are taken before the robots starts moving, while in the second  the channel estimation is updated
while the robot is in motion, and hence it is solved online. In both cases we utilize an in-house algorithm that computes near-optimal
solutions in little time, which enables its use in the online setting. The optimization strategy is described
in detail and validated by  simulation of realistic scenarios.
\end{abstract}


\section{Introduction}
{\it Communication-aware mobile robotics} is an emergent field of enquiry whose origins are in two related areas that
have been extensively researched in the past twenty years: mobile sensor networks \cite{cortes2002coverage,Mostofi_Chung_Murray_Burdick_IPSN05,leonard2007collective}, and networked robotics
systems \cite{BFPSS11,ZEP11,choi2009consensus,SSR02,marden2009cooperative}. A central problem in communication-aware robotics is in co-optimization of sensing, communication, and navigation
under physical and resource constraints. More specifically, the problem of balancing transmission energy
with motion energy has been the focus of research in recent years \cite{MMG10,MM11,GM11,YM10,YM14}, and it is also the subject of this paper.

While optimization of   transmission energy and motion energy
 has been traditionally explored separately in the respective literatures on communications and robotics (e.g.,    \cite{GC97,tokekar2013energy,MLHL06}), only recently
 the problem of co-optimizing the two forms of  energy has begun to attract attention. In \cite{OS09}, the authors propose an efficient approximate path planning algorithm that minimizes motion and communication energy costs. Ref. \cite{el2013mobile}
  optimizes  relay configurations in data-intensive wireless
 sensor networks.
  In \cite{tang2006energy}, the authors develop  an algorithm for maximizing
  the lifetime of  wireless sensor networks
 considering both  communication and motion costs of the sensors.   Ref.  \cite{jaleel2014Minimizing} considers a dynamic
 co-optimization problem in the setting of optimal control and develops a Hamiltonian-based algorithm for its solution.
  All of the references \cite{OS09,el2013mobile,tang2006energy,jaleel2014Minimizing} focus on the robotics
 and optimization aspects of the problem while using over-simplified models for channel and communication energy costs.
  Realistic models of channel fading \cite{MM11} are used in  \cite{yan2013co} and \cite{YM14} in designing
  a co-optimization strategy
 for balancing a robot's speed with transmission rate, and  \cite{usman2015optimal} develops an effective algorithm for realizing
 that strategy.

Ref. \cite{usman2015optimal} is the starting point of this paper. It considers a robot required to transmit a given
number of bits in a given amount of time to a remote station,  while traversing a predetermined path. The channel quality is variable
 along the path, and it is predicted by the realistic model described in  \cite{yan2013co}. The considered problem is
 to compute the profiles of the robot's acceleration and  spectral efficiency (transmission rate per unit bandwidth)
 that minimize the combined energy spent on transmission and motion. We cast the problem in
 the framework of optimal control, and solved it
 by using an in-house algorithm  that is simple to code and was shown  to yield fast convergence.

This paper extends the problem and methodology developed in \cite{usman2015optimal} in the following two ways. First, it adds the challenging element of path planning by lifting the restriction that the robot has to follow a pre-determined path. Rather, it has to
compute an optimal trajectory. Second, it considers an online scenario where the channel quality, generally obtained by measurements,
is updated while the robot is in motion thereby requiring a re-evaluation of the optimal trajectory
for the cost to go, in contrast to \cite{usman2015optimal} which does a one-time channel prediction
and co-optimization before the robot starts moving. We point out that we do not use model-predictive control or rolling horizons,
but rather compute the entire trajectory of the cost-to-go performance functional to the given
final time. These enhancements over the setting in \cite{usman2015optimal} pose significant computational challenges. The proposed scheme handles these challenges in an effective way, as shown by satisfactory solution for the example-problems presented in the sequel.

The rest of the paper is organized as follows. Section II  formulates the problem and discusses relevant existing results.
Section III, considering an elaborate example,  solves the  problem offline,  while section IV applies the online version of the algorithm.
Finally,  Section V concludes the paper and points out various directions for future research.
\vspace{0.1in}

\section{Problem Formulation} \label{PF}
This section formulates the optimal control problem, describes the channel prediction technique employed,
presents the algorithm used, and recounts results of its application to the fixed-path problem presented in
\cite{usman2015optimal}.

\subsection{Problem Definition.}\label{sec:prob_def}
Consider a robot  that has to traverse a path between
 a   source point $S\in{\mathcal R}^2$ and a  destination point $D\in{\mathcal R}^2$ while
transmitting a given number
of bits to a remote station in a given time-horizon $[0,t_f]$.
The problem is to determine the robot's path, acceleration, and transmission rate as
functions of time $t\in[0,t_f]$, that minimize the total energy required for transmission and motion.
  The power required for motion depends on the robot's  velocity and acceleration, while its transmission power
 depends on its position relative to the remote station, transmission rate,
  and the channel quality as will be detailed below.
  The channel quality is assessed by known estimation techniques that are based on spatial measurements
  which are performed by  various devices
  (such as static sensors in the field or crowdsourcing)  and transmitted to the robot periodically. Generally  the robot has to maintain a threshold reception-quality at the receiver (remote station) which requires a larger transmission power when the channel quality is lower.

The motion dynamics of the robot follow Newton's Law as follows,
\begin{align}	
	 & \dot{x}_1(t)= x_2(t),  \nonumber\\
	&\dot{x}_2(t)= u(t),
	\label{motion}
\end{align}
where $x_1 \in \mathcal{R}^2$ is the position of the robot in the plane, $x_2 \in \mathcal{R}^2$ denotes its velocity,
 and  $u \in \mathcal{R}^2$ is its  acceleration. The initial condition of this equation
 is $x_1(0)=S$ and $x_2(0)=0$.  According to Ref.   \cite{tokekar2013energy}, the power required for the robot's
 motion   has the form
\begin{align} 
P_m(t) &= k_1\vert \vert u(t)\vert \vert ^2 + k_2 \vert \vert x_2(t) \vert \vert^2 +k_3\vert \vert x_2(t)\vert \vert +k_4 \nonumber \\
	&+k_5\vert \vert u(t)\vert \vert +k_6\vert \vert u(t)\vert \vert \cdot \vert \vert x_2(t)\vert \vert,
\label{mobilitycost}
\end{align}
for given constants $k_i\geq 0$, $i=1, \ldots, 6$. The   power required for  transmitting data to the remote station is given by
\begin{align}
P_c (t) & = \frac{2^{R(t)}-1}{K} s(x_1(t)),
\label{eq:sx}
\end{align}
where $R(t)\geq 0$   is the spectral efficiency of the channel at time $t$ and position $x_1(t)$, $K$ is a constant depending on the threshold bit error rate acceptable at the receiver, and $s(x_1(t))$ is the estimated channel quality at position $x_1(t) \in \mathcal{R}^2$; see   Section \ref{sec:channel} for details. Let $Q$ be the total number of bits the robot has to transmit, and let $B$ denote the channel's bandwidth.
 Then the following equation represents the constraint that the robot has to transmit
 $Q$ bits during its motion,
\begin{align*}
	 \int_{0}^{t_f} R(t) dt  &= \frac{Q}{B}:=c.
\end{align*}
To represent this constraint without the integral, which would be more amenable to an application
 of the algorithm described below, we introduce an auxiliary  state variable, $x_3 \in \mathcal{R}$,
defined by the equation
\begin{align}
\dot{x}_3 &= R(t),
\label{data}
\end{align}
with the boundary conditions $x_3(0)=0$ and $x_3(t_f)=c$. Other final-time constraints
on the state variable (position and velocity) are $x_1(t_f)=D$ and $x_2(t_f)=0$.
We also assume  upper-bound constraints on $u(t)$ and $R(t)$ of the form
\begin{align}	
	&0 \leq \vert \vert u(t) \vert \vert  \leq u_{\max}, \nonumber\\
	&0\leq R(t)\leq R_{\max},
\label{input}	
\end{align}
for given $u_{\max}>0$ and $R_{\max}>0$.

The related  optimal control problem is defined as follows.
Its input is  $(u(t),R(t))\in{\mathcal R}^2\times{\mathcal R}$, $t\in[0,t_f]$, its state is $(x_1(t),x_2(t),x_3(t))$,
and its dynamics are given by Eqs. (\ref{motion}) and (\ref{data}) with the initial conditions
$x_1(0)=S$, $x_2(0)=0$, and $x_3(0)=0$. Its performance function, to be minimized, is
\begin{equation}
\bar{J}:= \int_0^{t_f}\big(P_{m}(t)+\gamma P_{c}(t)\big)dt,
\label{TotalCost}
\end{equation}
where $P_{m}(t)$ and $P_{c}(t)$ are the motion power and transmission power defined, respectively, by Eqs. (\ref{mobilitycost}) and (\ref{eq:sx}),
and $\gamma>0$ is a given constant. The
problem is to minimize $\bar{J}$ subject to the above dynamic equations, the upper-bound constraints on the input as defined
by Eq.  (\ref{input}), and the final-state constraints $x_1(t_f)=D$, $x_2(t_f)=0$,
and $x_3(t_f)=c$.

We handle the final-state constraints with a penalty function
of the form $C_1||x_1(t_f)-D||^2+C_2||x_2(t_f)||^2+C_3||x_3-c||^2$, for constants $C_1>0$, $C_2>0$, and $C_3>0$. The resulting optimal control problem
now has the following form: Minimize the cost functional $J$ defined as
	\begin{align}
	 J &= \int_{0}^{t_f} \bigg(\frac{2^{R(t)}-1}{K} s(x_1) + \gamma \big(k_1\vert \vert u(t)\vert \vert ^2 + k_2 \vert \vert x_2(t) \vert \vert ^2 \nonumber \\
	&+k_3\vert \vert x_2(t)\vert \vert +k_4+k_5\vert \vert u(t)\vert \vert +k_6\vert \vert u(t)\vert \vert\dot \vert \vert x_2(t)\vert \vert  \big)\bigg) dt \nonumber\\
	&+ C_1 \vert \vert x_1(t_f) - D \vert \vert ^2 +  C_2 \vert \vert x_2(t_f) \vert \vert ^2 + C_3 \vert \vert x_3(t_f) - c \vert \vert ^2,
	\label{cost2}
	\end{align}
subject to the dynamic equations
\begin{align*}
	& \dot{x}_1(t)= x_2(t), \hspace{0.3in} x_1(0)=S\\
	&\dot{x}_2(t)= u(t), \hspace{0.35in} x_2(0)=0\\
	& \dot{x}_3(t)= R(t), \hspace{0.35in} x_3(0)=0,
\end{align*}
and the constraints
\begin{align*}	
	& \vert \vert u(t) \vert \vert  \leq u_{\max}, \\
&	0\leq R(t)\leq R_{\max}.
\end{align*}

\subsection{Online Channel Prediction}\label{sec:channel}
Assuming the common MQAM modulation for a robot's communication to the remote station, the required transmit power at time $t$ can be characterized as \cite{G05}
\begin{align}
\tilde{P}_{\text{c}}(t)=(2^{R(t)}-1)/(K\Upsilon(x_1(t))),
\label{ActualcommPower}
\end{align}
where $K = - 1.5/\ln(5p_{b,\text{th}})$, $p_{b,\text{th}}$ is the given Bit Error Rate (BER) threshold at the receiver, $R(t)$ is the spectral efficiency at time $t$, $x_1(t) \in \mathcal{R}^2$ is the robot's position at time $t$, and $\Upsilon(x_1(t))$ the instantaneous channel-to-noise ratio (CNR) at $x_1(t)$. It is well known that the
CNR can be modeled as a random process with three components: path loss, shadowing and multipath fading\cite{G05}. As shown in \cite{MM11}, based on a small number of a priori channel measurements, a Gaussian random variable, $\Upsilon_{\text{dB}}(q)$, can best characterize the CNR (in the dB domain) at an unvisited location $q$, the mean and variance of which are given by
\begin{align*}
\overline{\Upsilon}_{\text{dB}}(q)= H_q \hat{\theta}+ \Psi^{\mathrm{T}}(q)\Phi^{-1}\big(Y -H_{\mathcal{Q}}\hat{\theta}), \\
\Sigma(q)= \hat{\xi}^2_{\text{dB}} + \hat{\rho}^2_{\text{dB}} -\Psi^{\mathrm{T}}(q)\Phi^{-1}\Psi(q),
\end{align*}
where $Y$ is the stacked vector of $m$ a priori gathered CNR measurements, $\mathcal{Q} =\{q_{1}, ... ,q_{m}\}$ denotes the measurement positions, $H_q = [1\;-10\log_{10}(\|q-q_{b}\|)]$, $H_{\mathcal{Q}} =[H_{q_1}^T\ ...\ H_{q_m}^T ]^T$, $\Phi = \Omega + \hat{\rho}_{\text{dB}}^{2}\:I_{m}$ with $[\Omega]_{i,j} = \hat{\xi}_{\text{dB}}^2 \exp(-\|q_{i} -q_{j}\|/\hat{\eta})$, for $i,j\in\{1, ... ,m\}$, and $\Psi(q)=[\hat{\xi}_{\text{dB}}^{2} \exp(-\|q-q_{1}\|/\hat{\eta})\ ...\ \hat{\xi}_{\text{dB}}^{2} \exp(-\|q -q_{m}\|/\hat{\eta})]^T$. The terms $\hat{\theta}=[\hat{K}_{\text{PL}} \ \hat{n}_{\text{PL}}]^T$, $\hat{\xi}_{\text{dB}}$, $\hat{\eta}$ and $\hat{\rho}_{\text{dB}}$ are the estimated channel parameters. See \cite{MM11} for more details on the estimation of channel parameters and the performance of this framework in channel prediction.

Based on this framework, the CNR at unvisited location $x_1(t)$ can be predicted as a lognormal random variable (in the linear domain). The expected trasmit power $P_{\text{c}}(t)$  is given by
\begin{align}
P_{\text{c}}(t) = \frac{2^{R(t)}-1}{K} E\left[\frac{1}{\Upsilon(x_1(t))}\right].
\label{commPower}
\end{align}
Note that for lognormally distributed $\Upsilon(x_1(t))$, we have
\begin{equation} \label{eq:comm_metric}
E\left[\frac{1}{\Upsilon(x_1(t))}\right] = \exp\left(\left(\frac{\ln 10}{10}\right)^2\frac{\Sigma(x_1(t))}{2}\right)\frac{1}{\overline{\Upsilon}(x_1(t))},
\end{equation}
where $\overline{\Upsilon}(x_1(t)) = 10^{\overline{\Upsilon}_{\text{dB}}(x_1(t))/10}$. Equation (\ref{eq:comm_metric}) provides an estimate of the predicted channel quality at $x_1(t)$ and we let $s(x_1(t)) = E\left[\frac{1}{\Upsilon(x_1(t))}\right]$, substituting which in Eq. (\ref{commPower}) leads to Eq. (\ref{eq:sx}) for computing the transmit power.

This framework is suitable to a setting where the channel prediction is updated as more channel measurements become available. Assume that the robot has a few channel measurement collected a priori (e.g. by static sensors in the field), based on which an initial prediction of channel quality over the workspace can be computed. The robot travels along the path obtained from minimizing
$J$ as defined in   (\ref{cost2}) with the initial channel prediction. As the robot moves, it is provided with additional  channel measurements (by gathering more samples along its path, through crowdsourcing and/or by other robots in
the field), which enables it to predict the
 channel quality more accurately. With such additional data, it  solves the problem  again for the remaining path,
 where in Eq. (\ref{cost2}) the starting time is the present time, say $t_{0}\in[0, t_f]$, and the initial condition (state) consists of the
 state $x(t_0)$ that has been obtained by the motion and transmission dynamics up to time $t_0$.  The details of this
  online optimization procedure will be presented   in Section \ref{sec:online opt}.

\subsection{Hamiltonian-Based Algorithm} \label{HBA}

Optimization algorithms typically are based on two computed   objects at a given iteration: a  direction (e.g., of descent), and a step size along
it. Recently we developed an algorithm which is suitable for a class of power-aware optimal control problems \cite{jaleel2014Minimizing}, \cite{usman2015optimal} and \cite{Hale14}.
Cumulative experience with it reveals some favorable computational properties including fast convergence towards
a local minimum. This does not mean fast asymptotic convergence, which characterizes an algorithm's
behavior close to a local minimum, but rather large strides towards  a region of a (local) minimum. A key innovation in
the algorithm is its choice of a descent direction, which is not based on gradient descent
but rather follows an alternative approach requiring little computing efforts.
We next explain the structure of the algorithm and  summarize its performance on
the power-aware problem considered in \cite{usman2015optimal}.

Consider the abstract Bolza optimal control problem where the system's dynamics
are defined by the equation
\begin{equation*}
\dot{x}=f(x,u)
\end{equation*}
with an initial condition $x(0):=x_0$,
where $x\in{\mathcal R}^n$, $u\in{\mathcal R}^k$, and $f:{\mathcal R}^n\times{\mathcal R}^k\rightarrow{\mathcal R}^n$
is Lipschitz continuous in $x$ and continuous in $u$.
Given a final time $t_f>0$, a running cost function $L:{\mathcal R}^n\times{\mathcal R}^k\rightarrow{\mathcal R}$, and
a terminal-state cost function $\phi:{\mathcal R}^n\rightarrow{\mathcal R}$, define the cost functional as
\begin{equation*}
J:=\int_{0}^{t_f}L(x,u)dt+\phi(x(t_f)).
\end{equation*}
The optimal control problem, considered, is to minimize $J$ subject to the pointwise constraints
$u(t)\in{\cal U}$, where ${\cal U}\subset{\mathcal R}^k$ is an input-constraint set. We make the following assumption:

\begin{assumption}
1). The function $f(x,u)$ is affine in $u\in{\cal U}$ for every $x\in{\mathcal R}^n$, and the
function $L(x,u)$ is convex in $u\in{\cal U}$ for every $x\in{\mathcal R}^n$. \\
2). The set ${\cal U}$ is compact and convex.
\end{assumption}

Let $p(t)$, $t\in[0,t_f]$, denote the costate (adjoint) trajectory defined by the equation
\begin{equation*}
\dot{p}=-\Big(\frac{\partial f}{\partial x}(x,u)\Big)^{\top}p-\Big(\frac{\partial L}{\partial x}\Big)^{\top}
\end{equation*}
with the boundary condition $p(t_f)=\nabla\phi(x(t_f)$, and let
\begin{equation*}
H(x,u,p):=p^{\top}f(x,u)+L(x,u)
\end{equation*}
denote the Hamiltonian function (see, e.g., Ref. \cite{bryson1975applied}).
The kind of problems for which our algorithm is suitable have the property that, for given $x\in{\mathcal R}^n$ and $p\in{\mathcal R}^n$,
a minimum value of the Hamiltonian  $H(x,w,p)$, over $w\in{\cal U}$, can be computed
via a simple, explicit formula.

 Any implementation of the algorithm would require approximations associated with computations
 performed only at time-points $t$ lying in a finite grid in the interval $[0,t_f]$. However,
 it is easier to describe the algorithm in its abstract, conceptual form where computations are performed at every $t\in[0,t_f]$. This is what
 we do in the present abstract description, but will specify implementation details, including the choice of a grid,
 when describing simulation experiments in the sequel. In the following description of the algorithm, we refer to the control function $u(t)$, $t\in[0,t_f]$, by the boldface notation
${\bf u}$.

\vspace{0.1in}

\noindent \textsc{Algorithm}

\vspace{0.05in}
{\it Parameters:} Constants $\alpha\in(0,1)$ and $\beta\in(0,1)$. Given a control ${\bf u}$, compute the next control ${\bf u}_{next}$ as follows:

\vspace{0.1in}
{\it 1. Direction from {\bf u}:}  Compute (numerically) the state and costate trajectories $x(t)$ and $p(t)$,
$t\in[0,t_f]$.
 For every $t\in[0,t_f]$, compute a pointwise ($t$-dependent) minimizer of the Hamiltonian denoted by
 $u^{\star}(t)$, namely,
a point $u^{\star}(t)\in{\cal U}$ satisfying
\begin{equation*}
u^{\star}(t)\in{\rm argmin}\Big(H(x(t),w,p(t))~|~w\in{\cal U}\Big).
\end{equation*}
Define ${\bf u}^{\star}$ to be the function $u^{\star}(t)$, $t\in[0,t_f]$.\footnote{There may arise measurability issues due to the
explicit characterization of $u^{\star}(t)$ for all $t$ in the uncountable set $[0,t_f]$. However, in a grid-based
implementation these issues will be avoided since $t$ would lie in a finite set.} Define the direction from ${\bf u}$ to
be $d(t):=u^{\star}(t)-u(t)$, namely, in functional notation,
$\textbf{d}=\textbf{u}^{\star}-\textbf{u}$.

\vspace{0.1in}
{\it 2. Step size along the direction ${\bf d}$:}
Define
\begin{equation*}
\theta(\textbf{u}) = \int_{0}^{t_f} \big(H(x(t),u^{\star}(t),p(t)) - H(x(t),u(t),p(t)) \big)dt.
\end{equation*}
Compute $k\in\{0,1,\ldots,\}$ defined as
\begin{equation}
k= \min \Big\{j=0,1, \ldots \big \vert \, J(\textbf{u} + \beta^j \textbf{d}) -J(\textbf{u}) \leq \alpha \beta^j \theta (\textbf{u})\Big\},
\label{Armijo}
\end{equation}
and set the step size, $\lambda$, to be $\lambda=\beta^k$.\\

{\it 3. Update:}
Set ${\bf u}_{next}$ to be
\begin{equation*}
\textbf{u}_{next} = \textbf{u} + \lambda\textbf{d}.
\label{unext}
\end{equation*}

As we pointed out, the main innovation of the algorithm is in the choice of the direction in Step 1,
while Step 2 describes a standard Armijo step size;  see \cite{polak1997optimization} for extensive discussions thereof.

The direction is not based on explicit gradient computations but rather comprises a form of conjugate gradient, and we
believe that this plays a role in the fast convergence of the algorithm that has been noted in various simulation
experiments. In particular, Fig. $2$ and Table I in \cite{usman2015optimal}, considering a one-dimensional
 power-aware problem, show that the algorithm makes most of its descent towards
the minimum value of $J$ in under  $10$ iterations. The next section exhibits similar
performance of the algorithm for the path-planning problem at hand.

\section{Path Planning with Motion and Communication Co-optimization} \label{PathPlanning}
In this section, we consider the application of the algorithm  to the problem defined in section II. The Hamiltonian associated with the  optimal control problem (\ref{cost2})  is
\begin{align}
H(x&,[u,R],p)= p_1^Tx_2+p_2^T u + p_3 R +  \frac{2^{R}-1}{K} s(x_1) +  \gamma \big(k_1\vert \vert u\vert \vert ^2  \nonumber\\
&\, \, \,+ k_2 \vert \vert x_2 \vert \vert ^2+k_3\vert \vert x_2\vert \vert +k_4+k_5\vert \vert u\vert \vert +k_6\vert \vert u\vert \vert \vert \vert x_2\vert \vert  \big),
\label{Hamiltonian}
\end{align}
where the costates $p_1 \in \mathcal{R}^2$ , $p_2 \in \mathcal{R}^2$, and $p_3 \in \mathcal{R}$ are defined by the
 adjoint equations
\begin{align*}
\dot{p}_1 &= - \frac{2^{R}-1}{K} \frac{\partial s(x_1)}{\partial x_1}, \nonumber \\
\dot{p}_2 &= -p_1 - \gamma \left(2k_2x_2 + k_3 \frac{x_2}{\vert \vert x_2\vert \vert} + k_6 \vert \vert u \vert \vert  \frac{x_2}{\vert \vert x_2\vert \vert} \right), \nonumber \\
\dot{p}_3 &= 0, 
\end{align*}
with terminal constraints $p_1(t_f)= 2C_1(x_1(t_f)-D)$, \quad $p_2(t_f)= 2C_2x_2(t_f)$ and $p_3(t_f)= 2C_3x_3(t_f)-c)$, respectively.

In the forthcoming we assume that $k_5=k_6=0$ in (\ref{Hamiltonian}) as we did in \cite{tokekar2013energy}. The minimizer of this Hamiltonian subject to the input constraints can be seen to be  given by
\begin{align}\nonumber
u^{\star} &= \begin{cases}
- \dfrac{p_2}{2 \gamma  k_1},  & \text{if } \dfrac{1}{2 \gamma k_1} \vert \vert p_2\vert \vert \leq u_{\max}\\
- \dfrac{p_2}{\vert \vert  p_2 \vert \vert} u_{max},   & \text{if } \dfrac{1}{2 \gamma k_1} \vert \vert p_2\vert \vert  > u_{\max},
\end{cases}
\end{align}

\begin{align}\nonumber
R^{\star} &= \begin{cases}
\dfrac{1}{\ln(2)} \ln\left(\dfrac{-p_3 K}{\ln(2) s(x_1)}\right), \, \, \text{if } p_3 \leq -\dfrac{(\ln(2) s(x_1))}{K}\\
R_{\max}, \hspace{0.8in}  \text{if } \dfrac{1}{\ln(2)} \ln\left(\dfrac{-p_3 K}{\ln(2) s(x_1)}\right) > R_{\max}\\
0, \hspace{1.2in} \text{otherwise}.\\
\end{cases}
\end{align}

\vspace{0.2in}
\noindent \textsc{Application:}

\vspace{0.1in}
Consider a robot that is tasked to move from the initial  point $S=(20,40)$ to the final point $D=(10,5)$ in the plane, and it has to transmit $150$ bits/Hz to a remote station located at $q_b=(5,5)$. The time budget available for the task is $40$ seconds. The acceleration and spectral efficiency can take maximum values of $u_{max}=0.5 m/s^2$ and $R_{max}=6 \, \frac{Bits/Hz}/{sec}$, respectively. The balancing factor between motion and communication was set to $\gamma=0.01$, and the constants $C_1$, $C_2$ and $C_3$ are set to $10$, $50$ and $10$, respectively. The Armijo step size parameters are set to $\alpha=0.1$ and $\beta=0.5$. The initial controls $u_0(t)$ and $R_0(t)$ are both set to zero. The integration  step size for the simulation is set to $dt=0.1$ seconds, and the algorithm is 
programmed to run for $500$ iterations. However, the  algorithm is terminated whenever the Armijo parameter $k$ in (\ref{Armijo}) is greater than $50$, indicating that a local minimum has been approached.

\begin{figure}[htbp]
\hspace{-0.1in}
\includegraphics[width=9cm, height=6.5cm]{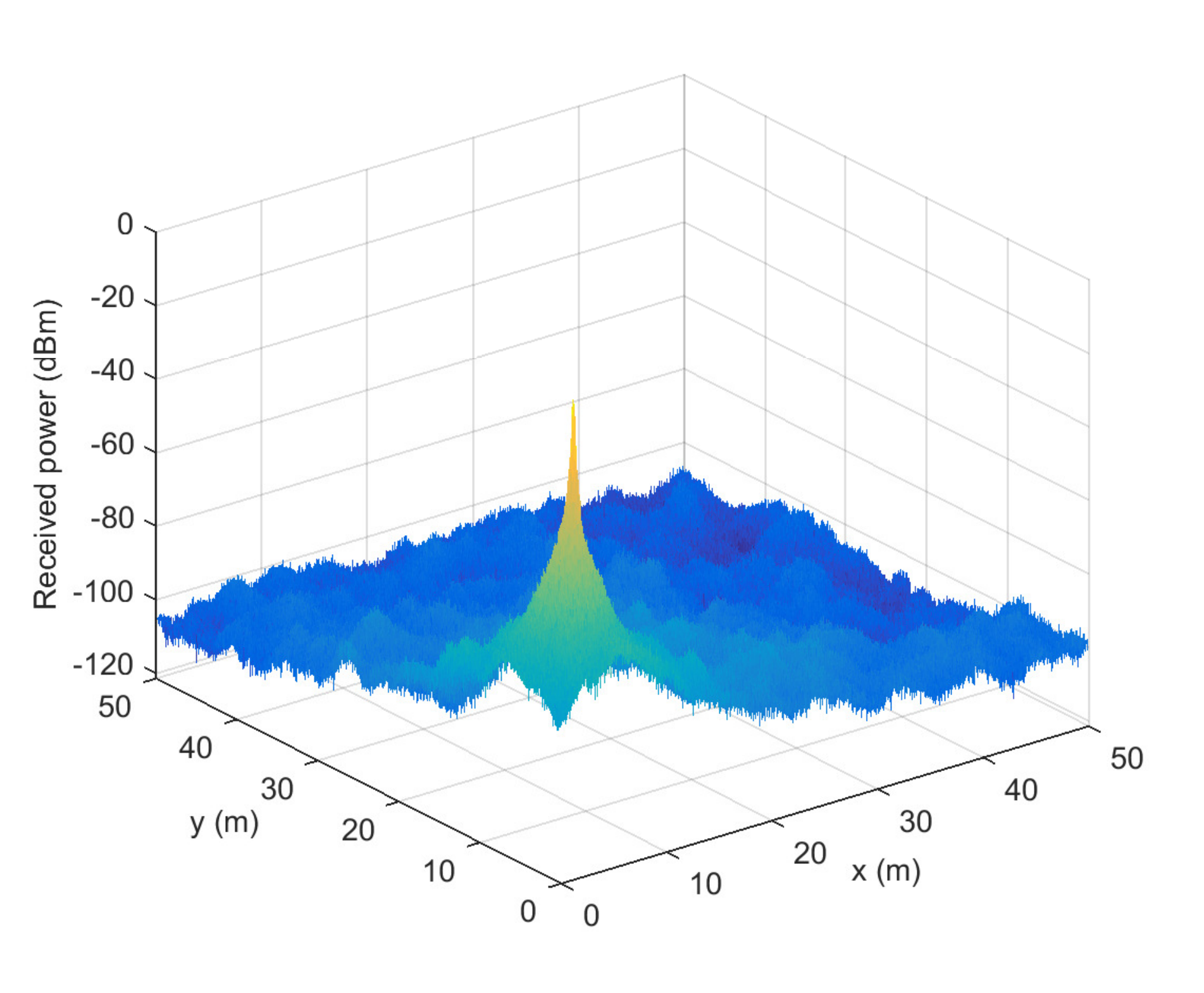}
\caption{Simulated wireless channel over the workspace. }
\label{figChannels1}
\end{figure}

\begin{figure}[htbp]
\hspace{-0.1in}
\includegraphics[width=9cm, height=6.5cm]{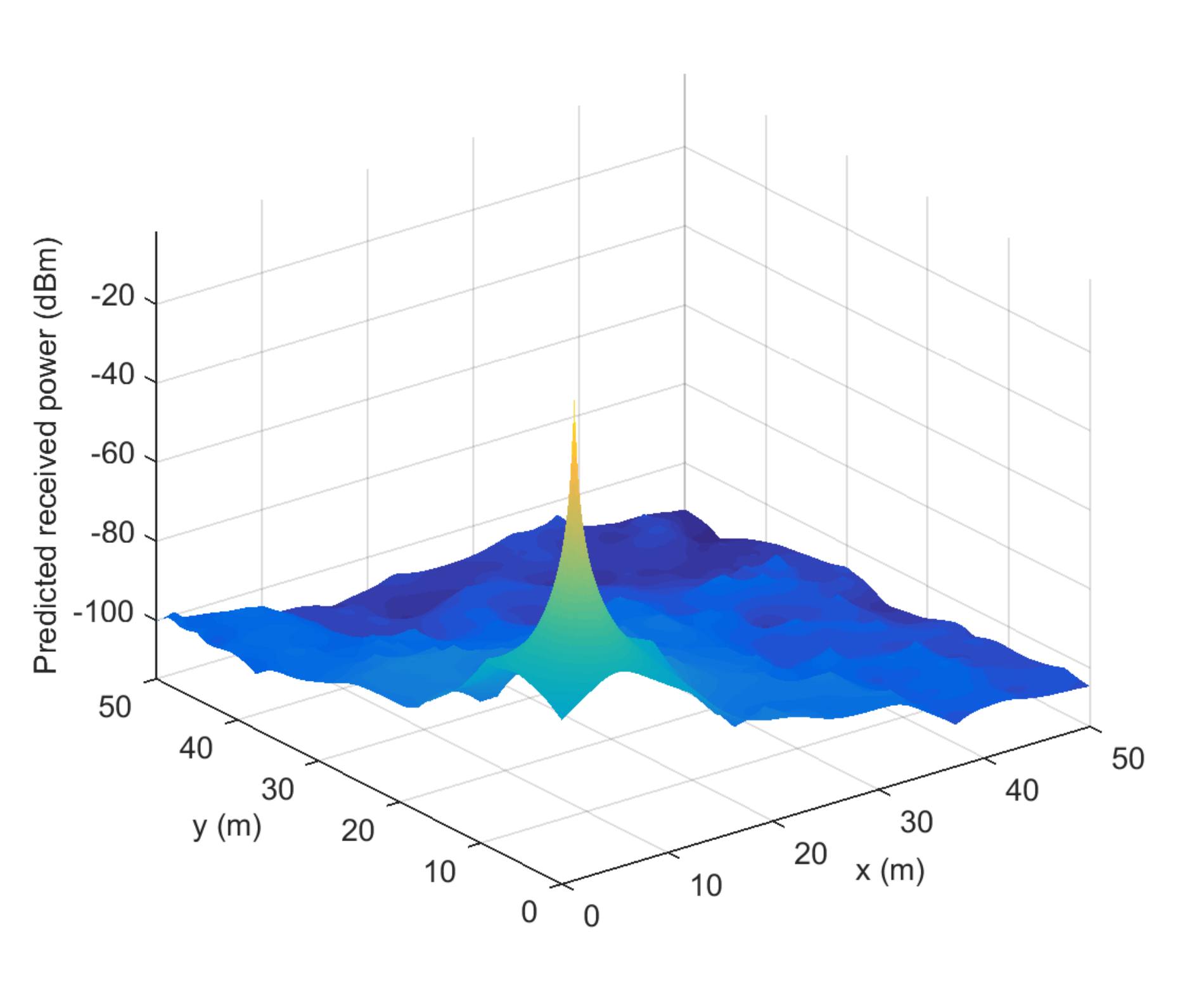}
\caption{Predicted channel based on $500$ measurements.}
\label{figChannels2}
\end{figure}

This robotic operation is performed under a simulated wireless channel with realistic parameters  over an area of $50m \times 50m$. The channel parameters based on \cite{MM11} and \cite{usman2015optimal}, are $K_{PL} = -41.34$, $n_{PL} = 3.86$, $\xi_{dB} = 3.20$, $\eta= 3.09 m$ and $\rho_{dB} = 1.64$. The receiver thermal noise is $-110$ dBm and the BER threshold is set to $p_{b,th}=2 \times 10^{-6}$. This channel can be predicted with  few  measurements over the field by using the methodology summarized in Section \ref{sec:channel}. To illustrate this point,  Fig. \ref{figChannels1} shows a sample simulated wireless channel generated with the parameters listed above, for the $250,000$ points in the plane. It is then predicted at all these points based on only $500$ a priori known randomly-spaced channel samples $(0.2\%)$ over the field and the result is shown in Fig. \ref{figChannels2}. The two results are quite similar.


\begin{figure}[htbp]
\hspace{-0.3in}
\includegraphics[width=10cm, height=7cm]{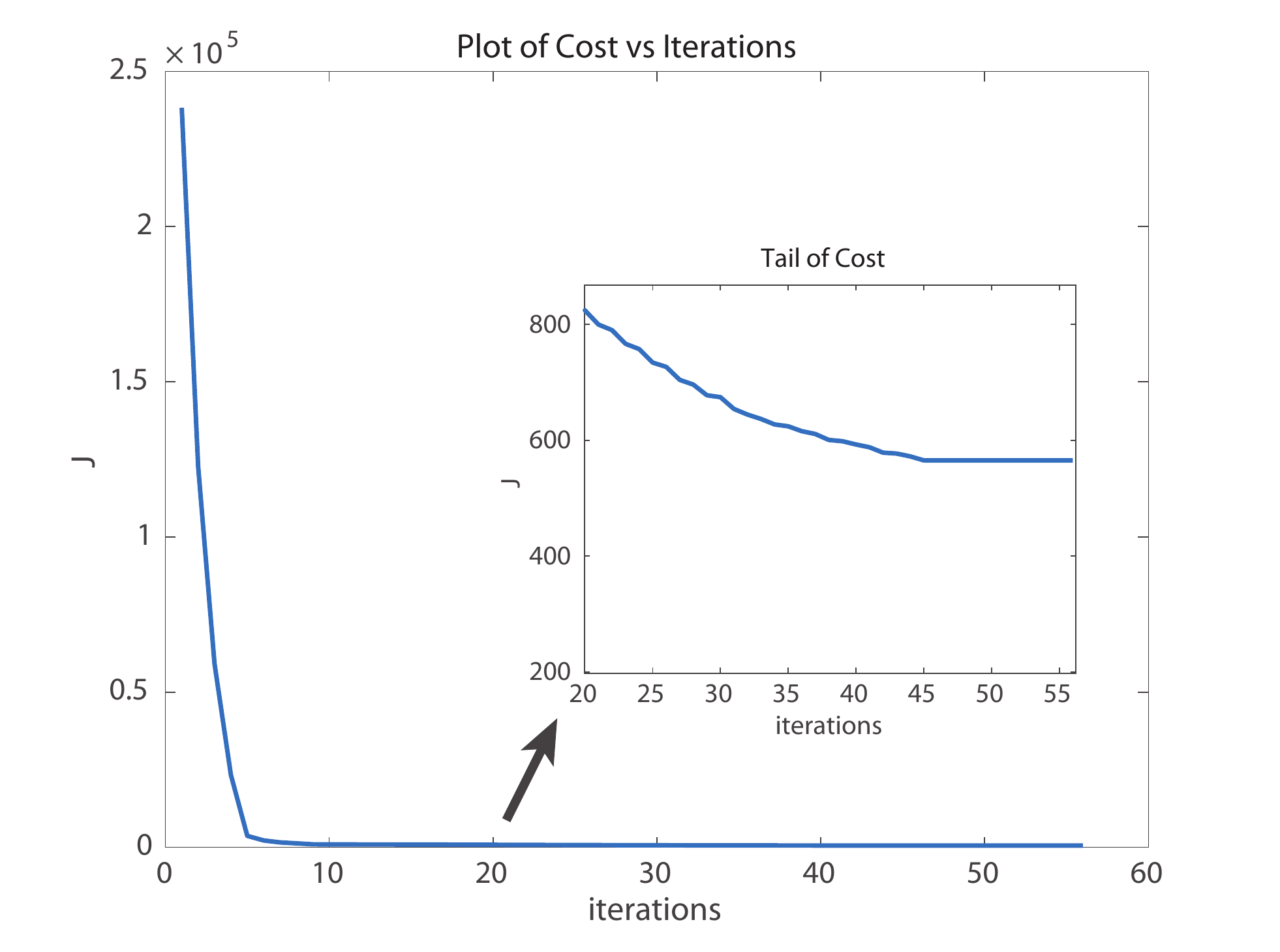}
\vspace{-0.3in}
\caption{Cost as function of iteration count.}
\label{J}
\end{figure}

A typical cost ($J$) vs. iteration count is depicted in Fig. \ref{J}.
Evidently the algorithm reaches values close to its obtained minimum in few iterations.
In fact, the computed cost is reduced  from the initial value of $2.3872 \times 10^{5}$ to $799.63$ after $20$ iterations,
while the cost after $56$ iterations is $565.13$  when $k$ became greater than $50$.  Fig. \ref{J} also  shows the tail of the cost trajectory and evidently it starts flattening after iteration 20. The $56$ steps of the algorithm took  $0.83$ seconds of CPU time on an Intel dual-core computer with i5 processor running at $2.7$ GHz.
%

The total motion and communication cost (\ref{TotalCost}), excluding the penalty term, is $\bar{J} = 475.10$, and 
the final values of state variables are $x_1(t_f)=(9.8, 5)$, $x_2(t_f) = (0.2,-0.8)$, and $x_3(t_f) =149.7$.  We note a
mild discrepancy from the  desired final values of $x_1(t_f)=(10,5)$, $x_2(t_f)=( 0,0)$, and $x_3(t_f)=150$,
 but it can be reduced by choosing larger penalty terms $C_1$, $C_2$, and $C_3$. As a matter of fact, a run of the algorithm with
   $C_1=500$, $C_2=500$, and $C_3=500$ gave  final states  of $x_1(t_f)=(9.99, 5)$, $x_2(t_f) = ( -0.08, -0.68)$, and $x_3(t_f) =149.99$; an initial cost of $J=1.1913 \times 10^7$,  and a final total mobility and communication cost of  $\bar{J} = 497.04$ after $500$ iterations. The  CPU time of the run was  $7.14$ seconds. It is not surprising that the 
   initial cost is higher since the penalty terms are larger, and for the same reason,  the
   algorithm drives the control parameters to a more restricted set and hence the final energy cost is expected to be higher
   as well. The CPU times often are larger in penalty-function methods with larger penalty terms.

\begin{figure}[htbp]
\begin{center}
\includegraphics[width=9cm, height=7.7cm]{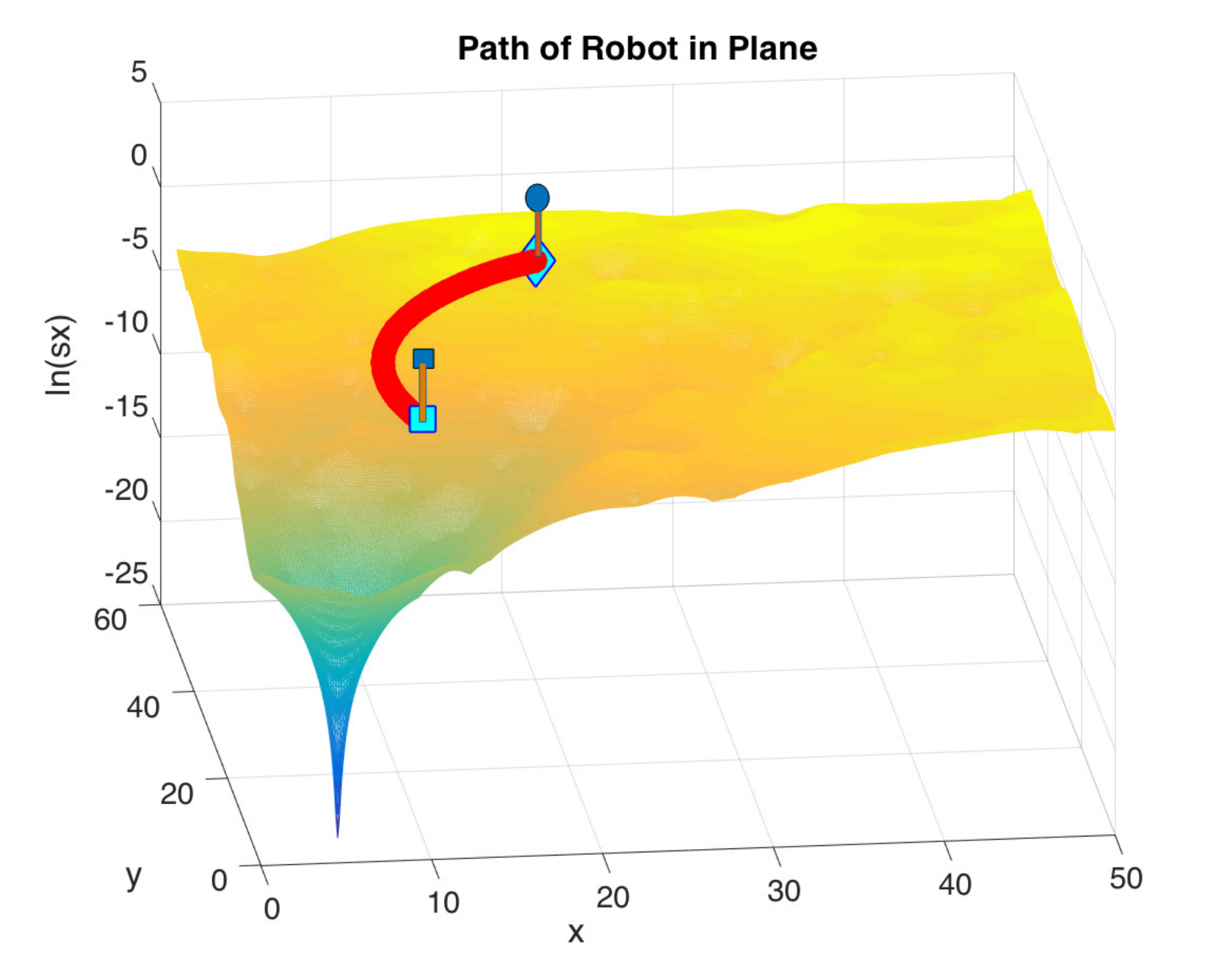}
\caption{Path followed by the robot, veering towards regions of better channel quality.}
\label{path}
\end{center}
\end{figure}

Fig. \ref{path} shows the log plot of predicted channel quality ($s(x_1)=E[1/\Upsilon(x,y)]$, where $\Upsilon(x,y)$ is the received CNR at position $x_1=(x,y)$) and the path  taken by the robot in the plane.  Smaller values of $s(x_1)$ correspond to good channel quality and vice versa.  The robot starting and end positions are marked as a square and a diamond,  respectively,  in all the figures. Instead of following a straight line between them, the robot takes a detour towards areas with predicted relatively good channel quality. For instance, the point of best channel quality is  $q_b=(5, 5)$, namely the location of base station, and hence the robot veers towards this point before turning away towards its destination point.

\begin{figure}[htbp]
\begin{center}
\includegraphics[width=9.2cm, height=7cm]{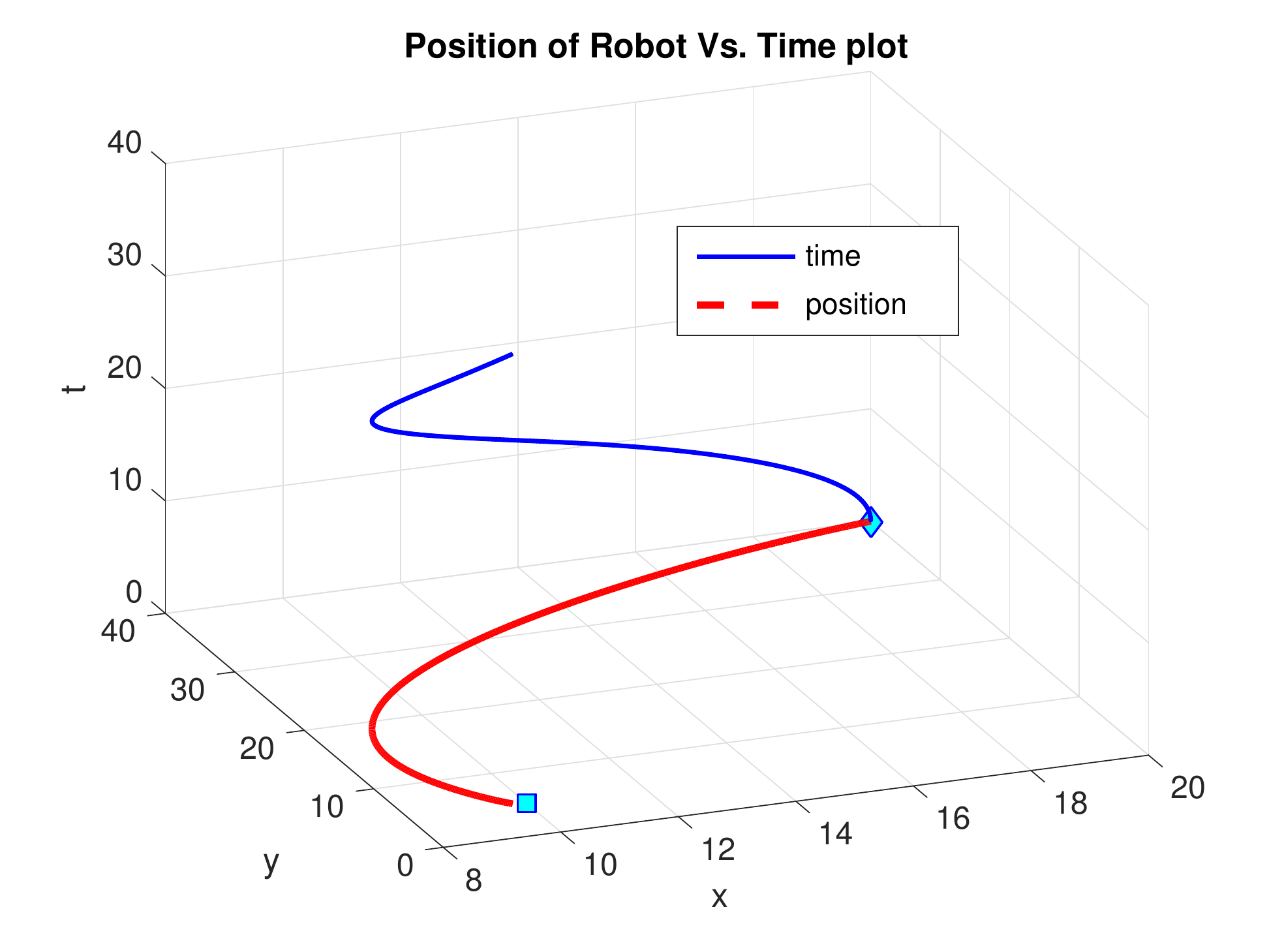}
\caption{Position of the robot as a function of time.}
\label{pathinplane}
\end{center}
\end{figure}

\begin{figure}[htbp]
\begin{center}
\includegraphics[scale=0.47]{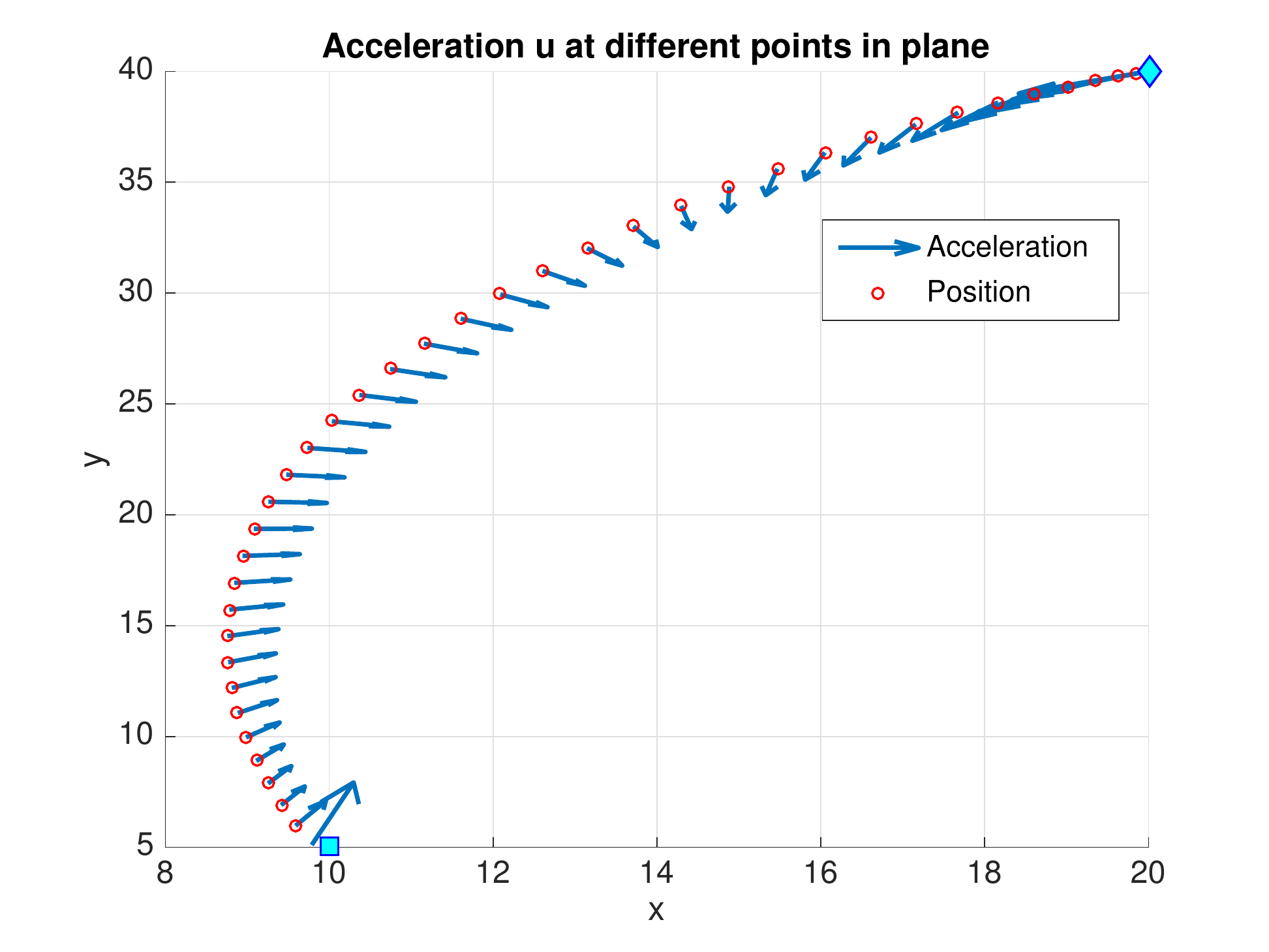}
\caption{Acceleration of the robot along its path.}
\label{acceleration}
\end{center}
\end{figure}

\begin{figure}[htbp]
\begin{center}
\includegraphics[width=9.2cm, height=7cm]{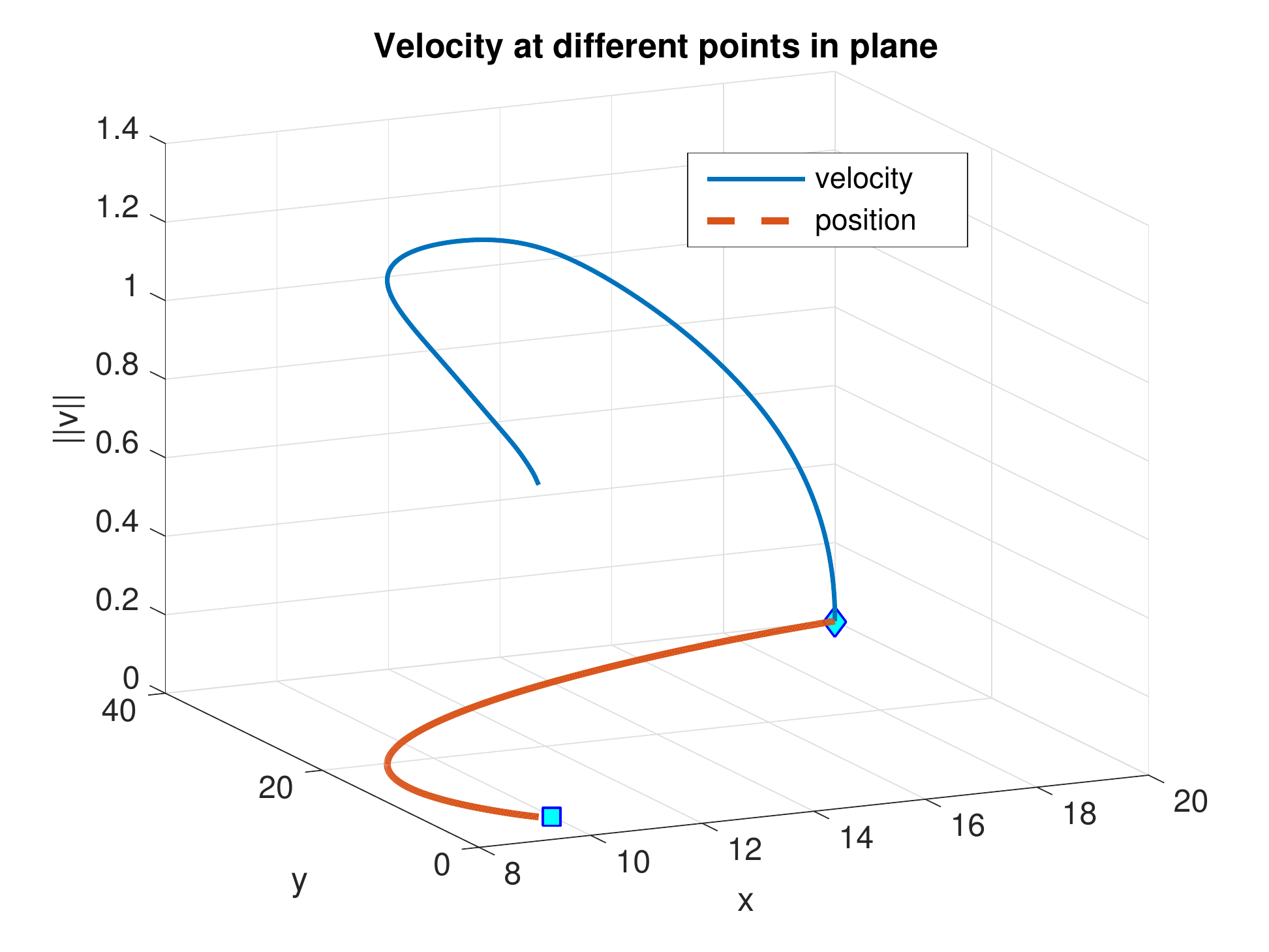}
\caption{Velocity of the robot along its path.}
\label{velocity}
\end{center}
\end{figure}

Fig. \ref{pathinplane} depicts a three-dimensional graph of the robot's motion, where the $t$ axis representing time and the 
 motion is in the $x-y$ plane.  The  upper, blue  curve represents the flow of time
 from $0$ to $40$ seconds, and the position of the robot at time $t$ is seen by projecting the corresponding point on the upper curve
 onto the $x-y$ plane, where it is indicated by a corresponding point
  on the red curve.   Fig. \ref{acceleration} shows the acceleration of the robot along its path, where   lengths of the arrows represent its magnitude, and  Fig. \ref{velocity} shows the speed of the robot along its path.

\begin{figure}[htbp]
\hspace{-0.15in}
\includegraphics[scale=0.47]{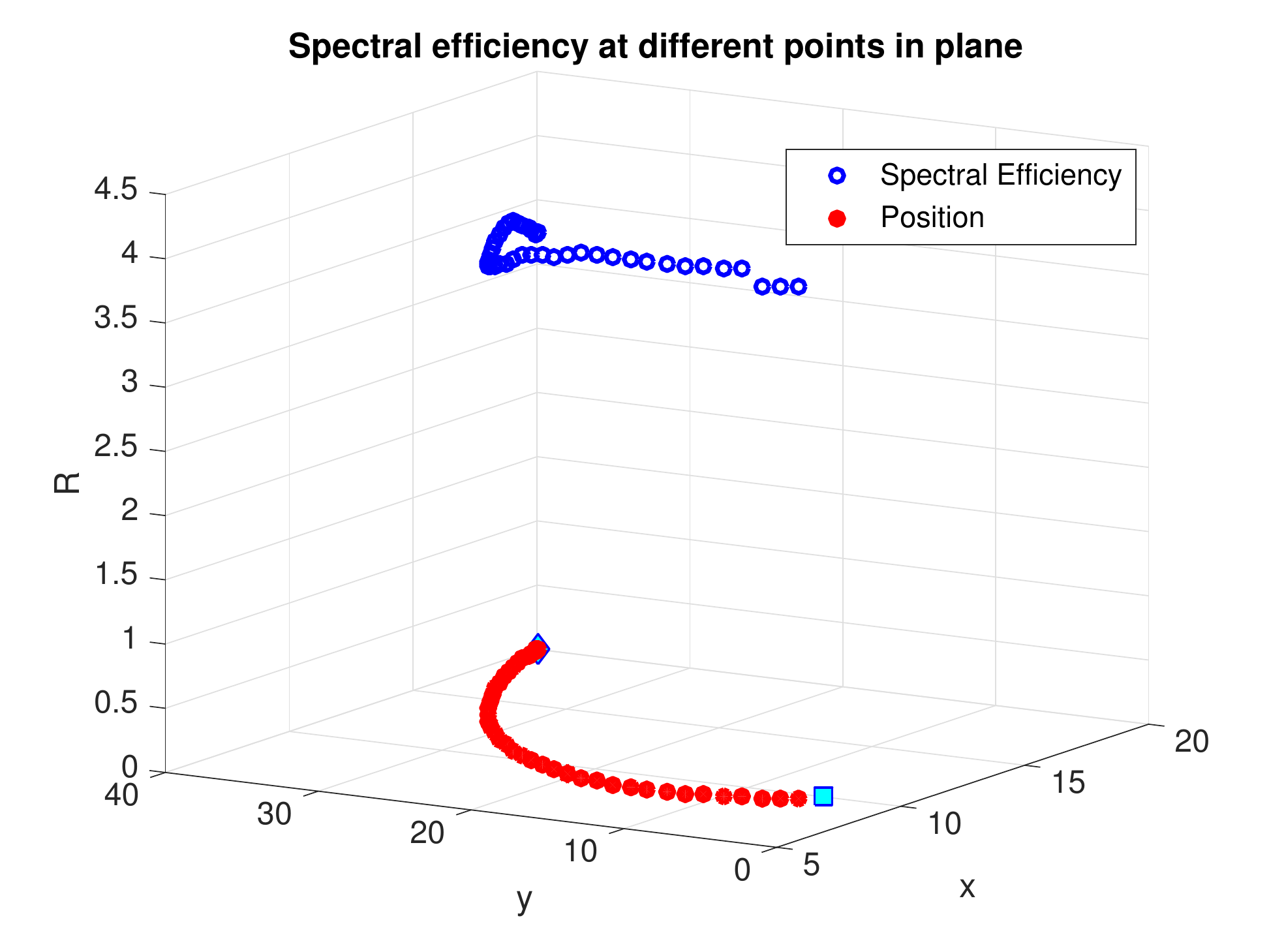}
\caption{The robot's spectral efficiency along its path.}
\label{Rx}
\end{figure}

\begin{figure}[htbp]
\hspace{-0.2in}
\includegraphics[width=9.2cm, height=6.5cm]{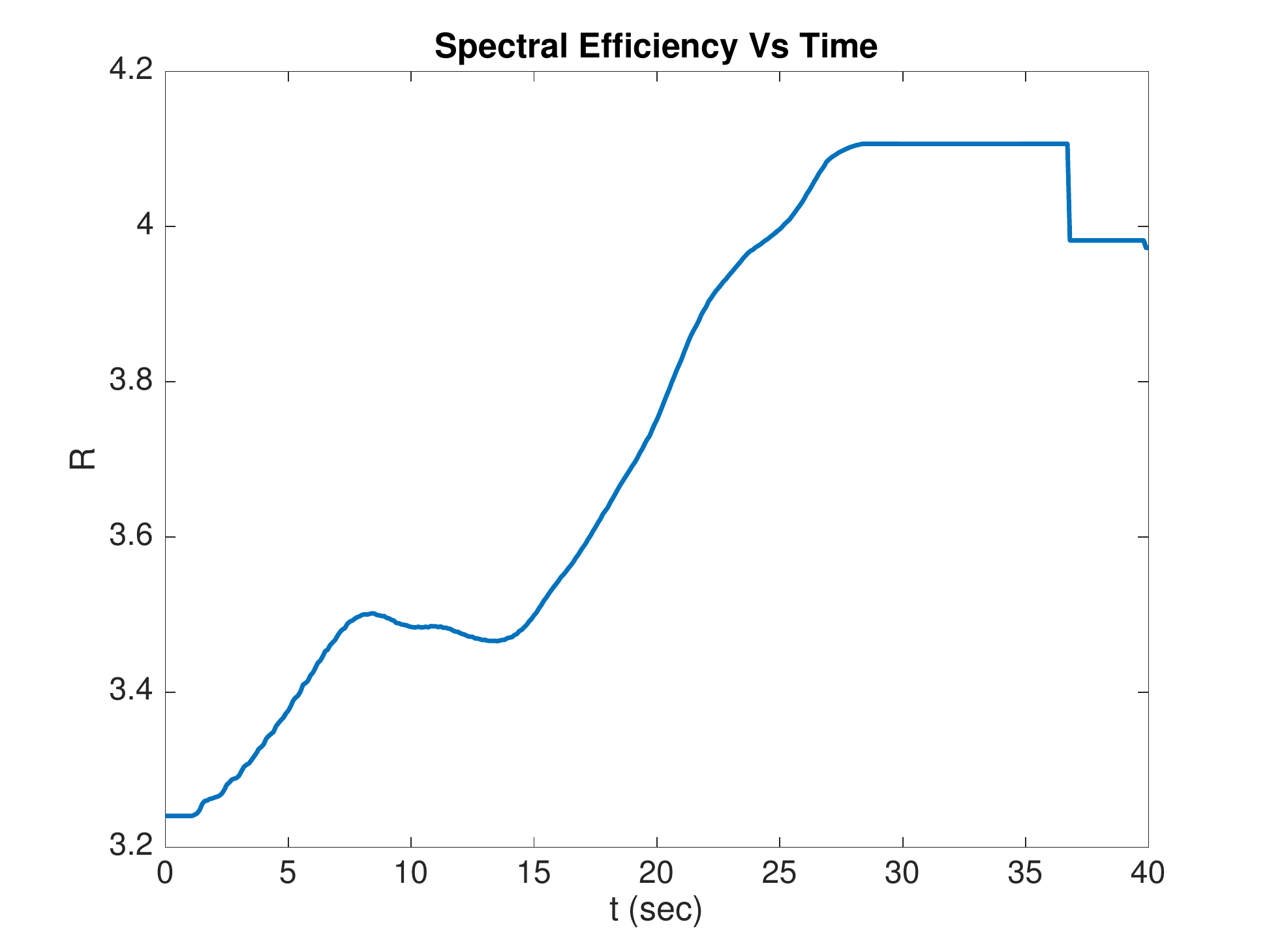}
\caption{The robot's spectral efficiency vs. time.}
\label{Rxt}
\end{figure}

The spectral efficiency of the robot is shown along its path in Figure  \ref{Rx}, where the path is marked by filled red circles,
  and the magnitude of spectral efficiency at corresponding points  is marked in blue. Fig. \ref{Rxt} shows the spectral efficiency vs time.  It can be seen that the robot transmits with higher spectral efficiency and hence at a  higher data rate in regions of better  channel quality.\footnote{Assuming a constant available  bandwidth, transmission rate is proportional to spectral efficiency.} This is not surprising since, in regions of higher channel quality, the robot can transmit more message bits  to the base station with less  communication power.
  
These figures show  that the robot veers towards regions of predicted
better channel-quality.  Within these regions  the required transmission power is  lower, and hence the robot slows down and increases spectral efficiency and  transmission rate in order  to send  a larger number of message bits.

\section{Online Optimization}\label{sec:online opt}
This section extends the algorithm to a realistic and practical  online setting, where the robot obtains new channel measurements while in
motion. It does not discard the older measurements, but rather appends them by the new data in order to enhance its
channel estimation.

The online optimal control problem is to minimize the cost functional $J_{t_0}$, defined as
	\begin{align}
	 &J_{t_0} = \int_{t_0}^{t_f} \bigg(\frac{2^{R(t)}-1}{K} s(x_1) + \gamma \big(k_1\vert \vert u(t)\vert \vert ^2 + k_2 \vert \vert x_2(t) \vert \vert ^2 \nonumber \\
	&+k_3\vert \vert x_2(t)\vert \vert +k_4+k_5\vert \vert u(t)\vert \vert +k_6\vert \vert u(t)\vert \vert \vert \vert x_2(t)\vert \vert  \big)\bigg) dt \nonumber\\
	&+ C_1 \vert \vert x_1(t_f) - D \vert \vert ^2 +  C_2 \vert \vert x_2(t_f) \vert \vert ^2 + C_3 \vert \vert x_3(t_f) - \bar{c} ) \vert \vert ^2,
	\label{cost3}
	\end{align}
subject to the dynamics 
\begin{align*}
	& \dot{x}_1(t)= x_2(t), \hspace{0.3in} x_1(t_0)=a_1\\
	&\dot{x}_2(t)= u(t), \hspace{0.35in} x_2(t_0)=a_2\\
	& \dot{x}_3(t)= R(t), \hspace{0.35in} x_3(t_0)=0
\end{align*}
and the constraints
\begin{align*}	
	&0 \leq \vert \vert u(t) \vert \vert  \leq u_{\max}, \\
	&0\leq R(t)\leq R_{\max},
\end{align*}
where $t_0\in[0,t_f]$ is the time at which the optimization is performed,   and the terms $a_1$ and $a_2$ are the current position and velocity of the robot at time $t_0$, and $\bar{c}:=(c -x_3(t_0^-))$  is the number of bits per unit frequency that remains to be transmitted in the time-interval $[t_0, \,t_f]$. The online approach solves this problem each time a channel estimation
is performed, typically at a finite number of times during the horizon $[0,t_f]$.  The initial control point of each such a run of the algorithm
consists of the remaining input control computed by its previous run. 

The considered problem is the one discussed in  Section \ref{PathPlanning}, except that  the robot performs channel
prediction  every $10$ seconds, and each prediction is based on $100$ new channel measurements taken at random locations.
Also the initial run, at $t_0=0$, solves the offline problem with $100$ channel samples. 
 The combined time for   channel prediction and a run of the algorithm  was about $2$ seconds and took under 50 iterations
 of the algorithm's run. 

\begin{figure}[htbp]
\hspace{-0.15in}
\includegraphics[scale=0.49]{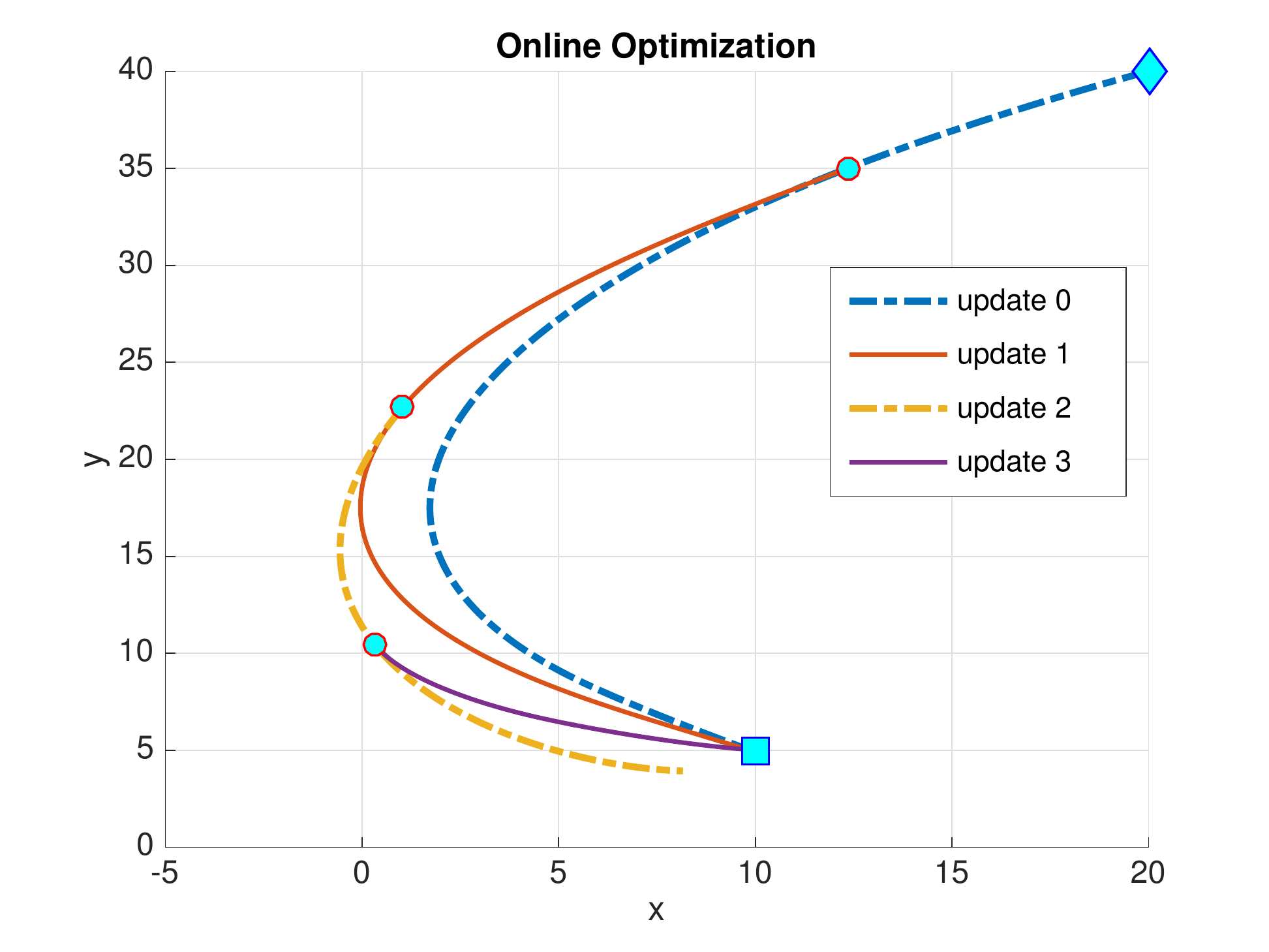}

\caption{Online optimization after every 10 seconds. }
\label{path1}
\end{figure}

\begin{figure}[htbp]
\hspace{-0.15in}
\includegraphics[scale=0.49]{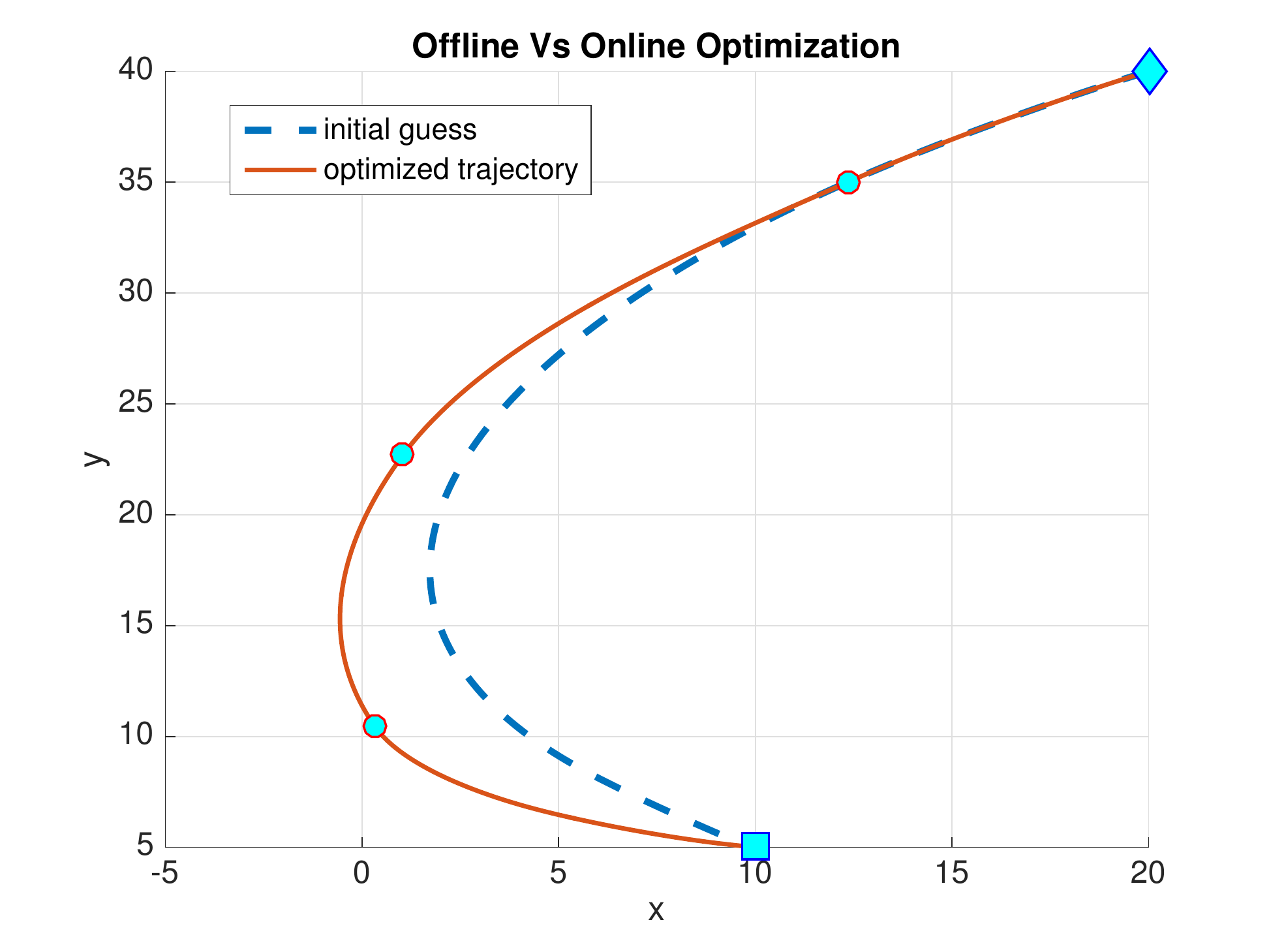}
\caption{Offline Vs. Online optimized Trajectories }
\label{path2}
\end{figure}

The results of the simulation are shown in Fig. \ref{path1} and Fig. \ref{path2}, where the position of the robot at the end of each predication and optimization cycle (10 seconds) is indicated by a circle. Fig. \ref{path1} shows the  computed optimal trajectories for each prediction-optimization cycle from the current time to the final time.  A concatenation of the computed trajectories,
 which the robot actually would traverse,  is indicated by the red  path in   Fig. \ref{path2}, while the dashed blue path 
  indicates the trajectory computed by the offline algorithm at time $t_0=0$, based on the initial channel prediction.
   The total energy consumed (Eq. (\ref{TotalCost})) in the offline solution  (dashed blue path Fig. \ref{path2}) is  $\bar{J}=371$,
    while the solution of the online problem (red path in Fig. \ref{path2}) yields a lower value,   $\bar{J}=304$.

\section{Conclusions and Future Research}
We considered the problem of co-optimization of communication and motion power of a robot that is required to transmit a given number of bits to a remote station in a pre-specified  amount of time, while moving between a given
starting point and an end point. The  problem is to compute
  the  robot's optimal path,and acceleration and transmission rate along this path.
   Both offline and online versions of the problem are considered and solved by simulation of realistic channel
   environments.
  Future work  would focus on extending the results from the case of a single robot
  to that  of  multiple agents having to perform coordinated tasks  while maintaining formation
  in the face of limited energy sources.

\vspace{0.2in}

\bibliographystyle{IEEEtran}
\bibliography{biblio-rev1}

\end{document}